\documentclass[12pt]{amsart}
\usepackage{mathrsfs}
\usepackage{cases}
\usepackage{epic}
\usepackage{amsfonts}
\usepackage{graphicx}
\usepackage{amsmath}
\usepackage{amssymb, upgreek}
\usepackage{bm}
\usepackage{latexsym,todonotes}
\usepackage{pdflscape}
\usepackage[all]{xypic}
\usepackage[all]{xy}
\usepackage{color}
\usepackage{colordvi}
\usepackage{multicol}
\usepackage[normalem]{ulem}
\usepackage[linktocpage=true]{hyperref}
\hypersetup{colorlinks,linkcolor=blue,urlcolor=cyan,citecolor=blue}

\begin{document}
\input xy
\xyoption{all}

\numberwithin{equation}{section}
\allowdisplaybreaks
\renewcommand{\mod}{\operatorname{mod}\nolimits}
\newcommand{\proj}{\operatorname{proj.}\nolimits}
\newcommand{\rad}{\operatorname{rad}\nolimits}
\newcommand{\soc}{\operatorname{soc}\nolimits}
\newcommand{\ind}{\operatorname{inj.dim}\nolimits}
\newcommand{\id}{\operatorname{id}\nolimits}
\newcommand{\Mod}{\operatorname{Mod}\nolimits}
\newcommand{\R}{\operatorname{R}\nolimits}
\newcommand{\End}{\operatorname{End}\nolimits}
\newcommand{\ob}{\operatorname{Ob}\nolimits}
\newcommand{\Ht}{\operatorname{Ht}\nolimits}
\newcommand{\cone}{\operatorname{cone}\nolimits}
\newcommand{\rep}{\operatorname{rep}\nolimits}
\newcommand{\Ext}{\operatorname{Ext}\nolimits}
\newcommand{\Tor}{\operatorname{Tor}\nolimits}
\newcommand{\Hom}{\operatorname{Hom}\nolimits}
\newcommand{\Pic}{\operatorname{Pic}\nolimits}
\newcommand{\aut}{\operatorname{Aut}\nolimits}
\newcommand{\Fac}{\operatorname{Fac}\nolimits}
\newcommand{\Div}{\operatorname{Div}\nolimits}
\newcommand{\rank}{\operatorname{rank}\nolimits}
\newcommand{\Len}{\operatorname{Length}\nolimits}
\newcommand{\RHom}{\operatorname{RHom}\nolimits}
\renewcommand{\deg}{\operatorname{deg}\nolimits}
\renewcommand{\Im}{\operatorname{Im}\nolimits}
\newcommand{\Ker}{\operatorname{ker}\nolimits}
\newcommand{\Iso}{\operatorname{Iso}\nolimits}
\newcommand{\Coh}{\operatorname{coh}\nolimits}
\newcommand{\Qcoh}{\operatorname{Qch}\nolimits}
\newcommand{\inj}{\operatorname{inj.dim}\nolimits}
\newcommand{\dimv}{\operatorname{\underline{\dim}}\nolimits}
\newcommand{\res}{\operatorname{res}\nolimits}

\def \U{\mathbf U}
\def \tUB{{}^{\texttt{B}}\tU}
\def \UB{{}^{\texttt{B}}\U}
\newcommand{\indim}{\operatorname{inj.dim}\nolimits}
\newcommand{\pdim}{\operatorname{proj.dim}\nolimits}
\def \dbl{\rm dbl}
\def \sqq{\mathbf v}
\newcommand{\Aut}{\operatorname{Aut}\nolimits}

\def \bm{\mathbf m}

\newcommand{\Cp}{\operatorname{Cp}\nolimits}
\newcommand{\coker}{\operatorname{Coker}\nolimits}
\renewcommand{\dim}{\operatorname{dim}\nolimits}
\renewcommand{\ker}{\operatorname{Ker}\nolimits}
\renewcommand{\div}{\operatorname{div}\nolimits}
\newcommand{\Ab}{{\operatorname{Ab}\nolimits}}
\newcommand{\Cone}{{\operatorname{Cone}\nolimits}}
\renewcommand{\Vec}{{\operatorname{Vec}\nolimits}}
\newcommand{\pd}{\operatorname{proj.dim}\nolimits}
\newcommand{\sdim}{\operatorname{sdim}\nolimits}
\newcommand{\add}{\operatorname{add}\nolimits}
\newcommand{\pr}{\operatorname{pr}\nolimits}
\newcommand{\oR}{\operatorname{R}\nolimits}
\newcommand{\oL}{\operatorname{L}\nolimits}

\newcommand{\tU}{\operatorname{\widetilde{\bf U}}\nolimits}

\newcommand{\Perf}{{\mathfrak Perf}}
\newcommand{\cc}{{\mathcal C}}
\newcommand{\ce}{{\mathcal E}}
\newcommand{\cs}{{\mathcal S}}
\newcommand{\cf}{{\mathcal F}}
\newcommand{\cx}{{\mathcal X}}
\newcommand{\cy}{{\mathcal Y}}
\newcommand{\cl}{{\mathcal L}}
\newcommand{\ct}{{\mathcal T}}
\newcommand{\cu}{{\mathcal U}}
\newcommand{\cm}{{\mathcal M}}
\newcommand{\cv}{{\mathcal V}}
\newcommand{\ch}{{\mathcal H}}
\newcommand{\ca}{{\mathcal A}}
\newcommand{\mcr}{{\mathcal R}}
\newcommand{\cb}{{\mathcal B}}
\newcommand{\ci}{{\mathcal I}}
\newcommand{\cj}{{\mathcal J}}
\newcommand{\cp}{{\mathcal P}}
\newcommand{\cg}{{\mathcal G}}
\newcommand{\cw}{{\mathcal W}}
\newcommand{\co}{{\mathcal O}}
\newcommand{\cd}{{\mathcal D}}
\newcommand{\ck}{{\mathcal K}}
\newcommand{\calr}{{\mathcal R}}

\def \fg{{\mathfrak g}}
\newcommand{\ol}{\overline}
\newcommand{\ul}{\underline}
\newcommand{\cz}{{\mathcal Z}}
\newcommand{\st}{[1]}
\newcommand{\ow}{\widetilde}
\renewcommand{\P}{\mathbf{P}}
\newcommand{\pic}{\operatorname{Pic}\nolimits}
\newcommand{\Spec}{\operatorname{Spec}\nolimits}
\newtheorem{theorem}{Theorem}[section]
\newtheorem{acknowledgement}[theorem]{Acknowledgement}
\newtheorem{algorithm}[theorem]{Algorithm}
\newtheorem{axiom}[theorem]{Axiom}
\newtheorem{case}[theorem]{Case}
\newtheorem{claim}[theorem]{Claim}
\newtheorem{conclusion}[theorem]{Conclusion}
\newtheorem{condition}[theorem]{Condition}
\newtheorem{conjecture}[theorem]{Conjecture}
\newtheorem{construction}[theorem]{Construction}
\newtheorem{corollary}[theorem]{Corollary}
\newtheorem{criterion}[theorem]{Criterion}
\newtheorem{propdef}[theorem]{Definition-Proposition}

\newtheorem{definition}[theorem]{Definition}
\newtheorem{example}[theorem]{Example}
\newtheorem{exercise}[theorem]{Exercise}
\newtheorem{lemma}[theorem]{Lemma}
\newtheorem{notation}[theorem]{Notation}
\newtheorem{problem}[theorem]{Problem}
\newtheorem{proposition}[theorem]{Proposition}
\newtheorem{solution}[theorem]{Solution}
\newtheorem{summary}[theorem]{Summary}
\newtheorem*{thm}{Theorem}
\newcommand{\qbinom}[2]{\begin{bmatrix} #1\\#2 \end{bmatrix} }

\theoremstyle{remark}
\newtheorem{remark}[theorem]{Remark}

\def \bfk{\mathbf k}
\def \bp{{\mathbf p}}
\def \bA{{\mathbf A}}
\def \bL{{\mathbf L}}
\def \bF{{\mathbf F}}
\def \bS{{\mathbf S}}
\def \bC{{\mathbf C}}
\def \bD{{\mathbf D}}
\def \Ire{\I^{\rm re}}
\def \Iim{\I^{\rm im}}
\def \Iiso{\I^{\rm iso}}

\def \I{\mathbb I}

\def \Z{{\Bbb Z}}
\def \F{{\Bbb F}}
\def \C{{\Bbb C}}
\def \N{{\Bbb N}}
\def \Q{{\Bbb Q}}
\def \G{{\Bbb G}}
\def \X{{\Bbb X}}
\def \P{{\Bbb P}}
\def \K{{\Bbb K}}
\def \E{{\Bbb E}}
\def \A{{\Bbb A}}
\def \BH{{\Bbb H}}
\def \T{{\Bbb T}}
\newcommand{\bluetext}[1]{\textcolor{blue}{#1}}
\newcommand{\redtext}[1]{\textcolor{red}{#1}}
\newcommand{\red}[1]{\redtext{ #1}}
\newcommand{\blue}[1]{\bluetext{ #1}}
\def \tMH{{\cs\cd\widetilde{\ch}}}

\title[Quantum Borcherds-Bozec algebras via semi-derived Ringel-Hall algebras]{Quantum Borcherds-Bozec algebras via semi-derived Ringel-Hall algebras}

\author[Ming Lu]{Ming Lu}
\address{Department of Mathematics, Sichuan University, Chengdu 610064, P.R.China}
\email{luming@scu.edu.cn}

	\dedicatory{Dedicated to Professor Liangang Peng on the occasion of his 65th birthday}

\subjclass[2020]{Primary 17B37, 
16E60, 18E35.}
\keywords{Quantum Borcherds-Bozec algebras, Semi-derived Ringel-Hall algebras, Quivers with loops}

\begin{abstract}
We use semi-derived Ringel-Hall algebras of quivers with loops to realize the whole quantum Borcherds-Bozec algebras and quantum generalized Kac-Moody algebras.
\end{abstract}

\maketitle
\section{Introduction}

In 1990, Ringel \cite{Rin90} constructed a Hall algebra associated to a Dynkin quiver $Q$ over a finite field $\mathbb F_q$, and identified its generic version with half a quantum group $\U^+ =\U^+_v(\fg)$; see Green \cite{Gr} for an extension to acyclic quivers. Ringel's construction has led to a geometric construction of $\U^+$ by Lusztig, who in addition constructed its canonical basis \cite{L90}. These constructions can be regarded as earliest examples of categorifications of halves of quantum groups. 

Bridgeland \cite{Br13} in 2013 used a Hall algebra of complexes to realize the whole quantum group $\U$. Actually Bridgleland's construction naturally produces the Drinfeld double $\widetilde \U$, a variant of $\U$ with the Cartan subalgebra doubled (with generators $K_i, K_i'$, for $i \in \I$). A  reduced version, which is the quotient of $\widetilde \U$ by the ideal generated by the central elements $K_i K_i'-1$, is then identified with $\U$. 

Bridgeland's version of Hall algebras has found further generalizations and improvements which allow more flexibilities. Gorsky \cite{Gor18} constructed {\em semi-derived Hall algebras} for Frobenius categories. More recently, motivated by the works of Bridgeland and Gorsky, the author and Peng \cite{LP16} formulated the ($\Z_2$-graded) {\em semi-derived Ringel-Hall algebras} starting with hereditary abelian categories.

We consider the quiver $Q$ with loops in this paper. Its Euler form yields a symmetric Borcherds-Cartan matrix $A$ and then a generalized Kac-Moody algebra $\fg_{A,\bm}$ with charge $\bm$ (also called Borcherds algebra) \cite{Bor88}. The quantum deformation of a generalized Kac-Moody algebra $\UB:=\UB_v(\fg_{A,\bm})$ and its modules are introduced in \cite{K95}. The Ringel-Hall algebra realization of the positive half $\UB^+$ of a quantum generalized Kac-Moody algebra is given in \cite{KS06}, where they also constructed the canonical basis theory.

Associated to a Borcherds-Cartan matrix, the Borcherds-Bozec algebra $\fg$ is a further generalization of the Kac-Moody algebra, which is
generated by infinitely many generators. 
The quantum Borcherds-Bozec algebra $\U:=\U_v(\fg)$ is introduced by Bozec in research of the perverse sheaves of a quiver with loops \cite{B15,B16}, which can be viewed as a further generalization of the quantum generalized Kac-Moody algebra. The Ringel-Hall algebra realization of the positive half $\U^+$ of a quantum Borcherds-Bozec algebra is given in \cite{K18} by using a quiver with loops.

In this paper, we use the semi-derived Ringel-Hall algebras of quivers with loops to realize the whole quantum Borcherds-Bozec algebras and quantum generalized Kac-Moody algebras, a generalization of Brigeland's construction. Similar to the quantum groups, our construction produces the Drinfeld double $\tU$ and $\tUB$, with their Cartan subalgebras doubled. The quotients of $\tU$ and $\tUB$ by the ideal generated by the central elements $K_i K_i'-1$ give the  quantum Borcherds-Bozec algebra $\U$ and quantum generalized Kac-Moody algebra $\UB$.

We give a remark here.
The definition of quantum Borcherds-Bozec algebras depends on parameters $\tau_{il}$ (see \eqref{eq:tau}) for any $(i,l)\in\I^\infty$. Only the quantum Borcherds-Bozec algebras with specific parameters $\tau_{il}$ \eqref{eq:tauspec} (up to isomorphisms) can be realized by semi-derived Ringel-Hall algebras.

\vspace{2mm}
\noindent{\bf Acknowledgments.} The author thanks Liangang Peng and Weiqiang Wang for guidance and continuing encouragement, and thanks Zhaobing Fan for helpful discussions on quantum Borcherds-Bozec algebras. 	The author also thanks University of Virginia, Shanghai Key Laboratory of Pure Mathematics and Mathematical Practice, East China Normal University for hospitality and support. This work is partially supported by the Science and Technology Commission of Shanghai Municipality (grant No. 18dz2271000), and the National Natural Science Foundation of China (grant No. 12171333).

\section{Semi-derived Ringel-Hall algebras}

Denote by $\Z_2:=\Z/2\Z$.
In this section, we shall review Ringel-Hall algebras, categories of $\Z_2$-graded complexes, semi-derived Ringel-Hall algebras for arbitrary hereditary abelian categories over a finite field $\bfk=\F_q$.

\subsection{Categories of $\Z_2$-graded complexes}
We assume that $\ca$  is an abelian category. 

Let $\cc_{\Z_2}(\ca)$ be the abelian category of $\Z_2$-graded complexes over $\ca$. Namely, an object $M$ of this category is a diagram with objects and morphisms in $\ca$:
$$\xymatrix{ M^0 \ar@<0.5ex>[r]^{d^0}& M^1 \ar@<0.5ex>[l]^{d^1}  },\quad d^1d^0=d^0d^1=0.$$
All indices of components of $\Z_2$-graded objects will be understood modulo $2$.
A morphism $s=(s^0,s^1):M\rightarrow N$ is a diagram
\[\xymatrix{  M^0 \ar@<0.5ex>[r]^{d^0} \ar[d]^{s^0}& M^1 \ar@<0.5ex>[l]^{d^1} \ar[d]^{s^1} \\
 N^0 \ar@<0.5ex>[r]^{e^0}& N^1 \ar@<0.5ex>[l]^{e^1} }  \]
with $s^{i+1}d^i=e^is^i$.

The shift functor on complexes is an involution
\begin{align}
\label{inv}
\xymatrix{\cc_{\Z_2}(\ca) \ar[r]^{*} & \cc_{\Z_2}(\ca)\ar[l],}
\end{align}
which shifts the grading and changes the sign of the differential as follows:
$$\xymatrix{ M^0 \ar@<0.5ex>[r]^{d^0}& M^1 \ar@<0.5ex>[l]^{d^1} \ar[r]^{*} & M^1 \ar@<0.5ex>[r]^{-d^1} \ar[l]& M^0 \ar@<0.5ex>[l]^{-d^0}  }.$$

Similar to ordinary complexes, we can define the $i$-th homology group for $M$, denoted by $H^i(M)$, for any $i\in\Z_2$. A complex is called acyclic if its homology group is zero. The subcategory formed by all acyclic complexes is denoted by $\cc_{\Z_2,ac}(\ca)$.

For any object $X\in\ca$, we define
\begin{align}
\label{stalks}
\begin{split}
K_X:=&(\xymatrix{ X \ar@<0.5ex>[r]^{1}& X \ar@<0.5ex>[l]^{0}  }),\qquad \,\, K_X^*:=(\xymatrix{ X \ar@<0.5ex>[r]^{0}& X \ar@<0.5ex>[l]^{1}  }),
\\
C_X:=&(\xymatrix{ 0 \ar@<0.5ex>[r]& X \ar@<0.5ex>[l]  }),\qquad \quad C_X^*:=(\xymatrix{ X\ar@<0.5ex>[r]& 0 \ar@<0.5ex>[l]  })
\end{split}
\end{align}
in $\cc_{\Z_2}(\ca)$. Note that $K_X,K_X^*$ are acyclic complexes.


From now on, we take the field $\bfk=\mathbb F_q$, a finite field of $q$ elements.
In the following, we always assume that $\ca$ is a hereditary abelian $\bfk$-linear category which is essentially small with finite-dimensional homomorphism and extension spaces.

\begin{lemma}[\text{\cite[Proposition 2.3]{LP16}}]
\label{proposition extension 2 zero}
Let $\ca$ be a hereditary abelian category.
For any $K\in\cc_{\Z_2,ac}(\ca)$,  we have
$$\pdim_{\cc_{\Z_2}(\ca)} (K)\leq 1\text{ and }\indim_{\cc_{\Z_2}(\ca)} (K)\leq 1.$$
\end{lemma}

Denote by $\Iso(\cc_{\Z_2}(\ca))$ the set of isomorphism classes $[M]$ of
$\cc_{\Z_2}(\ca)$.

By Proposition \ref{proposition extension 2 zero}, for any $[K],[M]\in \Iso(\cc_{\Z_2}(\ca))$ with $K$ acyclic, define
\begin{align*}
\langle [K],[M]\rangle=\dim\Hom_{\cc_{\Z_2}({\ca})}(K,M)-\dim\Ext^1_{\cc_{\Z_2}({\ca})}(K,M),
\\
\langle [M],[K]\rangle =\dim\Hom_{\cc_{\Z_2}({\ca})}(M,K)-\dim\Ext^1_{\cc_{\Z_2}({\ca})}(M,K).
\end{align*}
We call them the Euler forms. They descend to bilinear forms on the Grothendieck groups $K_0(\cc_{\Z_2,ac}(\ca))$
and $K_0(\cc_{\Z_2}(\ca))$, again denoted by the same symbol $\langle\cdot,\cdot\rangle$. 

Let $K_0(\ca)$ be the Grothendieck group of $\ca$.
For any $A\in\ca$, we denote by $\widehat{A}$ the corresponding element in  $K_0(\ca)$.
We also use $\langle \cdot,\cdot\rangle$ to denote the Euler form of $\ca$, i.e.,
\begin{align*}
\langle \widehat{A}, \widehat{B}\rangle=\dim\Hom_\ca(A,B)-\dim\Ext^1_\ca(A,B), \text{ for any }A,B\in\ca.
\end{align*}
Let $(\cdot,\cdot)$ be the symmetrized Euler form of $\ca$, i.e., $(\widehat{A},\widehat{B})= \langle \widehat{A}, \widehat{B}\rangle+\langle \widehat{B}, \widehat{A}\rangle$.


\begin{proposition}[\text{\cite[Proposition 2.4, Corollary 2.5]{LP16}}]
\label{lema euler form}
For any $A ,B \in\ca$, we have the following.
\begin{align}
&\langle [C_A],[K_B]\rangle=\langle [C_A^*],[K_B^*]\rangle=\langle \widehat{A}, \widehat{B}\rangle,\quad
\langle [K_B],[C_A^*]\rangle=\langle [K_B^*],[C_A]\rangle=\langle \widehat{B},\widehat{A}\rangle;
\\
&\langle [K_B],[C_A]\rangle=\langle [C_A^*],[K_B]\rangle=\langle [C_A],[K_B^*]\rangle=\langle [K_B^*],[C_A^*]\rangle=0;
\\
&\langle [K_{A}],[ K_{B}]\rangle= \langle [K_{A}^*], [K_{B}^*]\rangle=\langle [K_{A}], [K_{B}^*]\rangle=\langle [K_{A}^*], [K_{B}]\rangle=\langle \widehat{A},\widehat{B}\rangle.
\end{align}
\end{proposition}

Let $\ca\coprod\ca$ be the product of two copies of $\ca$. Then there is a forgetful functor $\res: \cc_{\Z_2}(\ca)\rightarrow \ca\coprod \ca$, which maps $M=(\xymatrix{ M^0 \ar@<0.5ex>[r]^{d^0}& M^1 \ar@<0.5ex>[l]^{d^1}  })$ to $(M^0, M^{1})$.



\subsection{Ringel-Hall algebras}

Let $\ca$ be an essentially small abelian category, linear over the finite field $\bfk=\F_q$.

Given objects $A,B,C\in\ca$, define $\Ext^1_\ca(A,C)_B\subseteq \Ext^1_\ca(A,C)$ to be the subset parameterising extensions with the middle term  isomorphic to $B$. We define the Hall algebra (also called Ringel-Hall algebra) $\ch(\ca)$ to be the $\Q$-vector space whose basis is formed by the isomorphism classes $[A]$ of objects $A$ of $\ca$, with the multiplication
defined by
\begin{align}
\label{eq:mult}
[A]\diamond [C]=\sum_{[B]\in \Iso(\ca)}\frac{|\Ext_\ca^1(A,C)_B|}{|\Hom_\ca(A,C)|}[B].
\end{align}
It is well known that
the algebra $\ch(\ca)$ is associative and unital. The unit is given by $[0]$, where $0$ is the zero object of $\ca$; see \cite{Rin90,Br13}. 

For any three objects $X,Y,Z$, let
\begin{align}
 \label{eq:Fxyz}
F_{XY}^Z:= \big |\{L\subseteq Z \mid L \cong Y,  Z/L\cong X\} \big |.
\end{align}
The Riedtmann-Peng formula states that
\[
F_{XY}^Z= \frac{|\Ext^1(X,Y)_Z|}{|\Hom(X,Y)|} \cdot \frac{|\Aut(Z)|}{|\Aut(X)| |\Aut(Y)|}.
\]

For any object $M$, 
let
\begin{align*}
[\![M]\!]:=\frac{[M]}{|\Aut(M)|}.
\end{align*}
Then the Hall multiplication \eqref{eq:mult} can be reformulated to be
\begin{align}
[\![M]\!]\diamond [\![N]\!]=\sum_{[\![L]\!]}F_{M,N}^L[\![L]\!],
\end{align}
which is the version of Hall multiplication used in \cite{Rin90}.

\subsection{Semi-derived Ringel-Hall algebras}

Let 
$$\sqq=\sqrt{q}.$$
Let $\widetilde{\ch}(\cc_{\Z_2}(\ca))$ be the twisted Ringel-Hall algebra of $\cc_{\Z_2}(\ca)$ over $\Q(\sqq)$, that is, $\widetilde{\ch}(\cc_{\Z_2}(\ca))$ has a basis formed by the isomorphism classes $[M]$ of objects $M$ of $\cc_{\Z_2}(\ca)$, with the product given by
\begin{align}\label{multiplication formula}
[L]* [M]=&\sqq^{\langle \res L,\res M\rangle}[L]\diamond [M]
\\\notag
=&\sqq^{\langle \res L,\res M\rangle}\sum_{[X]\in \Iso(\cc_{\Z_2}(\ca))}\frac{|\Ext^1_{\cc_{\Z_2}(\ca)}(L,M)_X|}{|\Hom_{\cc_{\Z_2}(\ca)}(L,M)|}[X].
\end{align}

Let $I_{\Z_2}$ be the two-sided ideal of $\widetilde{\ch}(\cc_{\Z_2}(\ca))$ generated by all differences $[L]-[K\oplus M]$ if there is a short exact sequence
\begin{equation}
  \label{eq:ideal}
 0 \longrightarrow K \longrightarrow L \longrightarrow M \longrightarrow 0
\end{equation}
with $K$ acyclic. 

Let $\widetilde{\ch}(\cc_{\Z_2}(\ca))/I_{\Z_2}$ be the quotient algebra. 
We also denote by $*$ the induced multiplication in $\widetilde{\ch}(\cc_{\Z_2}(\ca))/I_{\Z_2}$. In the following, we shall use the same symbols both in $\widetilde{\ch}(\cc_{\Z_2}(\ca))$ and $\widetilde{\ch}(\cc_{\Z_2}(\ca))/I_{\Z_2}$.
Let
 \begin{align}
 \label{multiset}
S_{\Z_2}:=\{a[K] \in \widetilde{\ch}(\cc_{\Z_2}(\ca))/I_{\Z_2} \mid  a \in \Q(\sqq)^\times, K\in\cc_{\Z_2,ac}(\ca)\},
 \end{align}
which is a multiplicatively closed subset with the identity $[0]\in S_{\Z_2}$.

\begin{propdef}[\cite{LP16}]
\label{proposition localizaition of Ringel-Hall algebra}
The multiplicatively closed subset $S_{\Z_2}$ is a right Ore, right reversible subset of $\widetilde{\ch}(\cc_{\Z_2}(\ca))/I_{\Z_2}$. Equivalently, the right localization of $\widetilde{\ch}(\cc_{\Z_2}(\ca))/I_{\Z_2}$ with respect to $S_{\Z_2}$ exists, which is called the twisted semi-derived Ringel-Hall algebra of $\ca$, and denoted by $\tMH(\ca)$.
\end{propdef}

Let $\ct(\ca)$ be the subalgebra of $\tMH(\ca)$ generated by $[K]^{\pm1}$ for all acyclic complexes $K$. Then $\ct(\ca)$ is a commutative algebra; see \cite[\S4.3]{LP16}.
Moreover, $\ct(\ca)$ is isomorphic to the group algebra of the Grothendieck group $K_0(\cc_{\Z_2,ac}(\ca))$.

For any $\alpha\in K_0(\ca)$, there exist $A,B\in\ca$ such that $\alpha=\widehat{A}-\widehat{B}\in K_0(\ca)$, and we set
\begin{align}
[K_\alpha]:=[K_A]* [K_B]^{-1},\qquad [K_\alpha^*]:=[K_A^*]*[K_B^*]^{-1}
\end{align}
Then $[K_\alpha],[K_\alpha^*]$ are well defined in $\tMH(\ca)$ (and then in $\ct(\ca)$); see \cite[\S3.4]{LP16}.

\begin{theorem}[cf. \text{\cite[Theorem 3.20]{LP16}}]
\label{theorem basis of modified hall algebra}
We have the following.
\begin{enumerate}
\item
$\tMH(\ca)$ has a basis given by
$$\{[C_A\oplus C_B^*]* [K_\alpha]*[K_\beta^*]\mid [A],[B]\in\Iso(\ca),\alpha,\beta\in K_0(\ca)\}.$$
\item
Let $M=(\xymatrix{ M^0 \ar@<0.5ex>[r]^{f^0}& M^1 \ar@<0.5ex>[l]^{f^1}  })\in \cc_{\Z_2}(\ca)$. Then we have, in  $\tMH(\ca)$,
\begin{align}
\label{complex to product}
[M]
=&\sqq^{ -\langle\widehat{H^0(M)}-\widehat{H^1(M)},\widehat{\Im f^0}-\widehat{\Im f^1}\rangle}
\, [C_{H^0(M)}^*\oplus C_{H^1(M)}]*[K_{\Im f^0}]* [K_{\Im f^1}^*].
\end{align}
\end{enumerate}
\end{theorem}

\begin{proof}
\cite[Theorem 3.20]{LP16} gives the corresponding results for the untwisted semi-derived Ringel-Hall algebras. The results here follow by considering the twisted Hall multiplication.
\end{proof}


Let $\widetilde{\ch}(\ca)$ be the twisted Ringel-Hall algebra, that is, the same vector space as $\ch(\ca)$ equipped with the twisted multiplication
$$[A]* [B]=\sqq^{\langle A,B\rangle}[A]\diamond [B]$$
for $[A],[B]\in\Iso(\ca)$.
Then the maps
\begin{align}
\label{eq:Rpm}
R^+:\widetilde{\ch}(\ca)&\longrightarrow \tMH(\ca),\qquad\qquad R^-:\widetilde{\ch}(\ca)\longrightarrow \tMH(\ca),
\\\notag
[A]&\mapsto [C_A],\qquad\qquad\qquad\qquad\qquad [A]\mapsto [C_A^*],
\end{align}
are algebra embeddings.
Moreover, by \cite[Theorem 3.25]{LP16}, we have the following triangular decomposition (as linear spaces)
\begin{align}
\label{eq:triangdecomp}
\tMH(\ca)\cong \widetilde{\ch}(\ca)\otimes \ct(\ca)\otimes \widetilde{\ch}(\ca).
\end{align}

\begin{lemma}
[\text{\cite[Lemma 4.3]{LP16}}]
\label{lem:KC}
For any $\alpha,\beta\in K_0(\ca)$, $M\in \ca$, we have
\begin{align}
\big[[K_\alpha],[K_\beta]\big]=&\big[[K_\alpha],[K_\beta^*]\big]=\big[[K_\alpha^*],[K_\beta^*]\big]=0,
\\
[K_\alpha] *[C_M]=&\sqq^{(\alpha,\widehat{M})}[C_M]* [K_\alpha], \qquad \qquad [K_\alpha] *[C_M^*]= \sqq^{-(\alpha,\widehat{M})} [C_M^*]*[K_\alpha],
\\
[K_\alpha^*]*[C_M]=&\sqq^{-(\widehat{M},\alpha)} [C_M]*[K_\alpha^*],\qquad\quad [K_\alpha^*]*[C_M^*]=\sqq^{(\widehat{M},\alpha)}[C_M^*]*[K_\alpha^*].
\end{align}
\end{lemma}

\section{Quantum Borcherds-Bozec algebras}
\label{sec:Drinfeld}

\subsection{Quantum Borcherds-Bozec algebras}
\label{subsec: QBB}
Let $\I$ be a finite index set. An (even) symmetric Borcherds-Cartan matrix is $A=(a_{ij})_{i,j\in \I}$ such that
\begin{enumerate}
\item $a_{ii}=2,0,-2,-4,\dots$,
\item $a_{ij}=a_{ji}\in\Z_{\leq0}$ for $i\neq j$.
\end{enumerate}
Denote by $\fg$ the Borcherds-Bozec algebra associated to $A$.

Let $\Ire:=\{i\in\I\mid a_{ii}=2\}$, $\Iim:=\{i\in\I\mid a_{ii}\leq0\}$, and $\Iiso:=\{i\in\I\mid a_{ii}=0\}$.
Let $\I^\infty:=(\Ire\times\{1\})\cup (\Iim\times \Z_{>0})$. For simplicity, we shall often write $i$ instead of $(i,1)$ for $i\in\Ire$.

Let $v$ be an indeterminant. Write $[A, B]=AB-BA$. Denote, for $r,m \in \N$,
\[
 [r]=\frac{v^r-v^{-r}}{v-v^{-1}},
 \quad
 [r]!=\prod_{i=1}^r [i], \quad \qbinom{m}{r} =\frac{[m][m-1]\ldots [m-r+1]}{[r]!}.
\]

For any $i\in\I$, denote by
\begin{align}
v_{(i)}=v^{\frac{a_{ii}}{2}}.
\end{align}

\begin{definition}
The algebra $\tU:=\tU_v(\fg)$ associated to a Borcherds-Cartan matrix $A$ is the associative algebra over $\Q(v)$ generated by the elements $K_i^{\pm1}$, $(K_i')^{\pm1}$ ($i\in\I$), $e_{il}$, $f_{il}$ ($(i,l)\in\I^\infty$) with the following defining relations:
\begin{align}
\label{eq:QBB1}
&K_iK_i^{-1}=K_i^{-1}K_i=1, \qquad K'_i(K'_i)^{-1}=(K'_i)^{-1}K'_i=1, \text{ for }i\in\I,
\\
\label{eq:QBB2}
&[K_i,K_j]=[K_i,K_j']=[K_i',K_j']=0, \text{ for }i,j\in\I,
\\
\label{eq:QBB3}
&K_ie_{jl}=v^{la_{ij}}e_{jl}K_i, \qquad
K_if_{jl}=v^{-la_{ij}}f_{jl}K_i,
\\
\label{eq:QBB4}
&K'_ie_{jl}=v^{-la_{ij}}e_{jl}K'_i, \qquad
K'_if_{jl}=v^{la_{ij}}f_{jl}K'_i,
\\
\label{eq:QBB5}
&e_{ik}f_{jl}-f_{jl}e_{ik}=0, \text{ for }i\neq j,
\\
\label{eq:QBB6}
&\sum_{\stackrel{m+r=k}{r+s=l}} v_{(i)}^{r(m-s)}\tau_{ir}e_{is}f_{im}(K_i')^r=\sum_{\stackrel{m+r=k}{r+s=l}} v_{(i)}^{-r(m-s)}\tau_{ir}f_{im}e_{is}(K_i)^r, \text{ for }(i,l),(i,k)\in\I^\infty,
\\
\label{eq:QBB7}
&
\sum_{k=0}^{1-la_{ij}}(-1)^k\qbinom{1-la_{ij}}{k} e_i^{1-la_{ij}-k}e_{jl}e_i^k=0,
\\
\label{eq:QBB8}
&\sum_{k=0}^{1-la_{ij}}(-1)^k\qbinom{1-la_{ij}}{k} f_i^{1-la_{ij}-k}f_{jl}f_i^k=0,
\text{ for }i\in\Ire, i\neq(j,l)\in\I^\infty,
\end{align}
where
\begin{align}
\label{eq:tau}
\tau_{il}\in 1+v\Z_{\geq0}[[v]],\qquad \forall (i,l)\in\I^\infty.
\end{align}
\end{definition}

For any $r\geq0$, define
\begin{align}
\varphi_r(t)=(1-t)(1-t^2)\cdots(1-t^r).
\end{align}
From now on, we shall assume that
\begin{align}
\label{eq:tauspec}
\tau_{ir}=\frac{1}{\varphi_r(v^2)}, \qquad\text{ for any }(i,r)\in\I^\infty.
\end{align}

Note that $K_iK_i'$ are central in $\tU$ for all $i\in\I$. Then the quantum Borcherds-Bozec algebra $\U:=\U_v(\fg)$ is defined to be the $\Q(v)$-algebra generated by the elements $K_i^{\pm1}$ ($i\in\I$), $e_{il}$, $f_{il}$ ($(i,l)\in\I^\infty$) subject to the defining relations \eqref{eq:QBB1}--\eqref{eq:QBB8} with $K_i'$ replaced by $K_i^{-1}$; see  \cite{B15,B16,FKKT20}.
In particular, $\U$ is isomorphic to the quotient algebra of $\tU$ modulo the ideal $(K_iK_i'-1\mid i\in\I)$.

Let $\tU^+$ be the subalgebra of $\tU$ generated by $e_{il}$ $((i,l)\in \I^\infty)$, $\tU^0$ be the subalgebra of $\tU$ generated by $K_i^{\pm1}, (K_i')^{\pm1}$ $(i\in \I)$, and $\tU^-$ be the subalgebra of $\tU$ generated by $f_{il}$ $((i,l)\in \I^\infty)$, respectively.
The subalgebras $\U^+$, $\U^0$ and $\U^-$ of $\U$ are defined similarly. Then both $\tU$ and $\U$ have triangular decompositions:
\begin{align}
\label{triandecomp2}
\tU =\tU^+\otimes \tU^0\otimes\tU^-,
\qquad
\U &=\U^+\otimes \U^0\otimes\U^-.
\end{align}
Clearly, ${\U}^+\cong\tU^+$, $\U^-\cong \tU^-$, and $\U^0 \cong \tU^0/(K_i K_i' -1 \mid   i\in \I)$.

\subsection{Semi-derived Ringel-Hall algebras of quivers}
Let $Q=(\I,\Omega)$ be a quiver, where $\I$ is the set of vertices and $\Omega$ is  the set of arrows. Define $g_i$ the number of loops at $i$ and $n_{ij}$ the number of arrows from $i$ to $j$ if $i\neq j$. Let $A=(a_{ij})_{i,j\in \I}$ be the symmetric Borcherds-Cartan matrix associated to $Q$, i.e.,
$$a_{ij}=\begin{cases}2-g_i& \text{ if }i=j, \\
-n_{ij}-n_{ji} & \text{ if }i\neq j.\end{cases}$$

A representation of $Q$ over $\bfk$ is $(M_i,x_\alpha)_{i\in\I,\alpha\in\Omega}$ by associating a finite-dimensional $\bfk$-linear space $M_i$ to each vertex $i\in\I$, and a  $\bfk$-linear maps $x_\alpha: M_i\rightarrow M_j$ to each arrow $\alpha:i\rightarrow j\in\Omega$.
A representation $(M_i,x_\alpha)$ of $Q$ is called nilpotent if $x_{\alpha_r}\cdots x_{\alpha_{2}}x_{\alpha_1}$ is nilpotent for any cyclic paths $\alpha_r\cdots \alpha_2\alpha_1$.

Let $\rep_\bfk^{\rm nil}(Q)$ be the category of finite dimensional nilpotent representations of $Q$ over $\bfk$. Then $\rep_\bfk^{\rm nil}(Q)$ is a hereditary abelian category.
Let $\langle\cdot,\cdot\rangle $ be the Euler form of $Q$. Define
$$(x,y)=\langle x,y\rangle+\langle y,x\rangle.$$
Let $S_i$ be the nilpotent simple module supported at $i\in\I$.
Then $(S_i,S_j)=a_{ij}$ for any $i,j\in\I$.

Denote by $\sqq_{(i)}=\sqq^{\frac{a_{ii}}{2}}=\sqq^{1-g_i}$ for $i\in\I$.
\begin{lemma}
For any $k,l\geq1$, in $\tMH(\rep_\bfk^{\rm nil}(Q))$ we have
\begin{align}
\label{eq:SkSl1}
[\![C_{S_i^{\oplus k}}]\!]*[\![C^*_{S_i^{\oplus l}}]\!]
=&\sum_{r=0}^{\min(k,l)} \frac{q^{-r(l+k-r)}}{\varphi_r(q^{-1})} \sqq^{r(l-k)}_{(i)} [\![C_{S_i^{\oplus (k-r)}}\oplus C^*_{S_i^{\oplus (l-r)}}]\!]*[K_{S_i}^*]^{ r},
\\
\label{eq:SkSl2}
[\![C^*_{S_i^{\oplus l}}]\!]*[\![C_{S_i^{\oplus k}}]\!]
=&\sum_{r=0}^{\min(k,l)} \frac{q^{-r(l+k-r)}}{\varphi_r(q^{-1})} \sqq^{r(l-k)}_{(i)} [\![C_{S_i^{\oplus (k-r)}}\oplus C^*_{S_i^{\oplus (l-r)}}]\!]*[K_{S_i}]^{ r}.
\end{align}
\end{lemma}

\begin{proof}
We only prove \eqref{eq:SkSl1} since the proof of \eqref{eq:SkSl2} is similar.

For any $k,l\geq1$, a direct computation shows
\begin{align*}
&[C_{S_i^{\oplus k}}]*[C^*_{S_i^{\oplus l}}]
\\
=&\sum_{r=0}^{\min(k,l)} \prod_{t=0}^{r-1}\frac{(q^k-q^t)(q^l-q^t)}{q^r-q^t} [C_{S_i^{\oplus (k-r)}}\oplus C^*_{S_i^{\oplus (l-r)}}\oplus (K_{S_i}^*)^{\oplus r}]
\\
=&\sum_{r=0}^{\min(k,l)}\sqq^{r(l+k-r)+\binom{r}{2}}(\sqq-\sqq^{-1})^r \qbinom{k}{r}_\sqq \qbinom{l}{r}_\sqq [r]_\sqq^![C_{S_i^{\oplus (k-r)}}\oplus C^*_{S_i^{\oplus (l-r)}}\oplus (K_{S_i}^*)^{\oplus r}]
\\
=&\sum_{r=0}^{\min(k,l)}\sqq^{r(l+k-r)+\binom{r}{2}}\sqq^{r(l-k)(1-g_i)} (\sqq-\sqq^{-1})^r \qbinom{k}{r}_\sqq \qbinom{l}{r}_\sqq [r]_\sqq^![C_{S_i^{\oplus (k-r)}}\oplus C^*_{S_i^{\oplus (l-r)}}]*[K_{S_i}^*]^{ r}.
\end{align*}
Here the second equality follows from
\begin{align*}
\prod_{t=0}^{r-1} (q^r -q^t) &=\sqq^{r^2 +{r \choose 2}} (\sqq -\sqq^{-1})^r [r]_\sqq^{!},
\\
\prod_{t=0}^{r-1} (q^k -q^t) &=\sqq^{rk +{r \choose 2}} (\sqq -\sqq^{-1})^r [k]_\sqq [k-1]_\sqq \ldots [k-r+1]_\sqq,
\end{align*}
and the last equality follows from \eqref{complex to product}.

Note that
\begin{align*}
&|\Aut(C_{S_i^{\oplus k}})|=\prod_{t=0}^{k-1} (q^k -q^t),\qquad
|\Aut(C^*_{S_i^{\oplus l}})|=\prod_{t=0}^{l-1} (q^l -q^t),
\\
&|\Aut(C_{S_i^{\oplus k}}\oplus C^*_{S_i^{\oplus l}})|=\prod_{t=0}^{k-1} (q^k -q^t)\cdot\prod_{t=0}^{l-1} (q^l -q^t).
\end{align*}
It follows that
\begin{align*}
&[\![C_{S_i^{\oplus k}}]\!]*[\![C^*_{S_i^{\oplus l}}]\!]
\\
=&\sum_{r=0}^{\min(k,l)}\sqq^{r(l+k-r)+\binom{r}{2}-k^2-\binom{k}{2}-l^2-\binom{l}{2}+(l-r)^2+\binom{l-r}{2}+(k-r)^2+\binom{k-r}{2}}\sqq^{r(l-k)(1-g_i)} (\sqq-\sqq^{-1})^{-r} \\
&\qquad\times\frac{1}{[r]_\sqq^!}[\![C_{S_i^{\oplus (k-r)}}\oplus C^*_{S_i^{\oplus (l-r)}}]\!]*[K_{S_i}^*]^{ r}
\\
=&\sum_{r=0}^{\min(k,l)}\sqq^{\frac{5r^2+r}{2}-2rk-2rl}\sqq^{r(l-k)(1-g_i)} (\sqq-\sqq^{-1})^{-r} \frac{1}{[r]_\sqq^!}[\![C_{S_i^{\oplus (k-r)}}\oplus C^*_{S_i^{\oplus (l-r)}}]\!]*[K_{S_i}^*]^{ r}
\\
=&\sum_{r=0}^{\min(k,l)} \frac{q^{-r(l+k-r)}}{\varphi_r(q^{-1})} \sqq^{r(l-k)}_{(i)} [\![C_{S_i^{\oplus (k-r)}}\oplus C^*_{S_i^{\oplus (l-r)}}]\!]*[K_{S_i}^*]^{ r}.
\end{align*}
Here we use $\sqq_{(i)}=\sqq^{1-g_i}$.
\end{proof}

\begin{lemma}
For any $k,l\geq1$, we have
\begin{align}
\label{eq:commha}
\sum_{r=0}^{\min(k,l)} \frac{q^{-r(k+l)+r(r+1)}}{\varphi_r(q)} \sqq_{(i)}^{r(l-k)}[\![C_{S_i^{\oplus (k-r)}}]\!]*[\![C^*_{S_i^{\oplus (l-r)}}]\!]*[K_{S_i}^*]^r
\\
\notag
=\sum_{r=0}^{\min(k,l)} \frac{q^{-r(k+l)+r(r+1)}}{\varphi_r(q)} \sqq_{(i)}^{r(l-k)}[\![C_{S_i^{\oplus (l-r)}}]\!]*[\![C^*_{S_i^{\oplus (k-r)}}]\!]*[K_{S_i}]^r.
\end{align}
\end{lemma}

\begin{proof}
We shall prove that
\begin{align}
\label{eq:commha1}
[\![C_{S_i^{\oplus k}}\oplus C^*_{S_i^{\oplus l}}]\!]=\sum_{r=0}^{\min(k,l)} (-1)^r\frac{q^{-r(k+l)+\binom{r+1}{2}}}{\varphi_r(q^{-1})} \sqq_{(i)}^{r(l-k)}[\![C_{S_i^{\oplus (k-r)}}]\!]*[\![C^*_{S_i^{\oplus (l-r)}}]\!]*[K_{S_i}^*]^r,
\\
[\![C_{S_i^{\oplus k}}\oplus C^*_{S_i^{\oplus l}}]\!]=\sum_{r=0}^{\min(k,l)} (-1)^r\frac{q^{-r(k+l)+\binom{r+1}{2}}}{\varphi_r(q^{-1})} \sqq_{(i)}^{r(l-k)}[\![C^*_{S_i^{\oplus (l-r)}}]\!]*[\![C_{S_i^{\oplus (k-r)}}]\!]*[K_{S_i}]^r,
\end{align}
and then \eqref{eq:commha} follows.

We only prove \eqref{eq:commha1} by induction on $k$, since the second one is similar.
For convenience, we assume $[\![C_{S_i^{\oplus k}}\oplus C^*_{S_i^{\oplus l}}]\!]=0$ if either $k$ or $l$ is negative. Similarly, $[\![C_{S_i^{\oplus k}}]\!]=0=[\![ C^*_{S_i^{\oplus k}}]\!]$ if $k<0$. It is direct to prove \eqref{eq:commha1} for $k=1$ by using \eqref{eq:SkSl1}. It follows by induction that
\begin{align*}
&[\![C_{S_i^{\oplus k}}\oplus C^*_{S_i^{\oplus l}}]\!]
\\
=&[\![C_{S_i^{\oplus k}}]\!]*[\![C^*_{S_i^{\oplus l}}]\!]-\sum_{r\geq1} \frac{q^{-r(l+k-r)}}{\varphi_r(q^{-1})} \sqq^{r(l-k)}_{(i)} [\![C_{S_i^{\oplus (k-r)}}\oplus C^*_{S_i^{\oplus (l-r)}}]\!]*[K_{S_i}^*]^{ r}
\\
=&[\![C_{S_i^{\oplus k}}]\!]*[\![C^*_{S_i^{\oplus l}}]\!]-\sum_{r\geq1}\sum_{p\geq0} \frac{q^{-r(l+k-r)}}{\varphi_r(q^{-1})} \sqq^{(r+p)(l-k)}_{(i)} (-1)^p\frac{q^{-p(k+l-2r)+\binom{p+1}{2}}}{\varphi_p(q^{-1})} \\
&\qquad\times[\![C_{S_i^{\oplus (k-r-p)}}]\!]*[\![C^*_{S_i^{\oplus (l-r-p)}}]\!]*[K_{S_i}^*]^{r+p}
\\
=&[\![C_{S_i^{\oplus k}}]\!]*[\![C^*_{S_i^{\oplus l}}]\!]-\sum_{a\geq1}(-1)^a  \sqq^{a(l-k)}_{(i)} q^{-a(l+k)+\binom{a+1}{2}} \Big( \sum_{r=1}^a (-1)^r q^{-\binom{r+1}{2}+ar} \frac{1}{\varphi_r(q^{-1})\varphi_{a-r}(q^{-1})}  \Big)
\\
&\qquad\times[\![C_{S_i^{\oplus (k-a)}}]\!]*[\![C^*_{S_i^{\oplus (l-a)}}]\!]*[K_{S_i}^*]^{a}
\\
=&\sum_{a=0}^{\min(k,l)} (-1)^a\frac{q^{-a(k+l)+\binom{a+1}{2}}}{\varphi_a(q^{-1})} \sqq_{(i)}^{a(l-k)}[\![C_{S_i^{\oplus (k-a)}}]\!]*[\![C^*_{S_i^{\oplus (l-a)}}]\!]*[K_{S_i}^*]^a.
\end{align*}
The last equality used the following formula (cf. \cite[1.3.1(c)]{Lus93}; see also \cite[Proof of Lemma 4.14]{LRW21})
$$\sum_{r=1}^a (-1)^r \frac{q^{-\binom{r+1}{2} +ar}}{\varphi_r(q^{-1})\varphi_{a-r}(q^{-1})} =-\frac{1}{\varphi_a(q^{-1})}.$$
\end{proof}

\begin{theorem}
There exists an injective morphism of algebras $\Phi: \tU\rightarrow \tMH(\rep^{\rm nil}_\bfk(Q))$ such that
\begin{align}
\Phi(K_i)=\sqq^2[K_{S_i}],\quad \Phi(K_i')=\sqq^2[K_{S_i}^*],\quad \Phi(e_{il})=\sqq^{l^2}[\![C_{S_i^{\oplus l}}]\!],\quad \Phi(f_{il})=\sqq^{l^2}[\![C^*_{S_i^{\oplus l}}]\!]
\end{align}
for any $i\in\I$, $(i,l)\in\I^\infty$.
\end{theorem}

\begin{proof}
In order to prove that the well-definedness  of $\Phi$, it is enough to check that $\Phi$ preserves the relations \eqref{eq:QBB1}--\eqref{eq:QBB8}. \eqref{eq:QBB1}--\eqref{eq:QBB4} follow from Lemma \ref{lem:KC}; \eqref{eq:QBB6} follows from \eqref{eq:commha};
for \eqref{eq:QBB5}, it follows by
$$[C_{{S_i}^{\oplus k}}]*[C^*_{{S_j}^{\oplus l}}]=[C_{{S_i}^{\oplus k}}\oplus C^*_{{S_j}^{\oplus l}}]=[C^*_{{S_j}^{\oplus l}}]*[C_{{S_i}^{\oplus k}}],\quad\text{ if }i\neq j;$$
\eqref{eq:QBB7}--\eqref{eq:QBB8} follow from \cite[Theorem 4.2]{K18} by using \eqref{eq:Rpm}.

For the injectivity of $\Phi$, note that
$\tU =\tU^+\otimes \tU^0\otimes\tU^-$ by \eqref{triandecomp2} and $\tMH(\rep^{\rm nil}_\bfk(Q))\cong \widetilde{\ch}(\rep^{\rm nil}_\bfk(Q))\otimes \ct(\rep^{\rm nil}_\bfk(Q))\otimes \widetilde{\ch}(\rep^{\rm nil}_\bfk(Q))$ by \eqref{eq:triangdecomp}. Moreover, we have
$\ct(\rep^{\rm nil}_\bfk(Q))\cong \Q(\sqq)\big[[K_{S_i}],[K_{S_i}^*]\mid i\in\I\big]$ by using $K_0(\rep_\bfk^{\rm}(Q))\cong\Z^\I$. Then the result follows from \cite[Theorem 4.2]{K18} by using the same argument of \cite[Theorem 4.9]{Br13}.
\end{proof}

\section{Quantum generalized Kac-Moody algebras}

In this section, we use semi-derived Ringel-Hall algebras of quivers with loops to realize the quantum generalized Kac-Moody algebras briefly, which strengthens \cite[Remark 4.16]{LP16}.

Let $A=(a_{ij})_{i,j\in\I}$ be a symmetric Borcherds-Cartan matrix as defined in \S\ref{subsec: QBB}. Assume that we are given a collection of positive integers $\bm=(m_i)_{i\in\I}$ with $m_i=1$ whenever $i\in\Ire$, called the charge of $A$. Let $\fg_{A,\bm}$ be the generalized Kac-Moody algebra (also called Borcherds algebra).

Let $\tUB:=\tUB_v(\fg_{A,\bm})$ be the $\Q(v)$-algebra generated by the elements $K_i^{\pm1}$, $(K_i')^{\pm1}$, $E_{ik}$, $F_{ik}$ for $i\in\I$, $k=1,\dots,m_i$ subject to \eqref{eq:QBB1}--\eqref{eq:QBB2} and
\begin{align}
\label{eq:QB3}
&K_iE_{jl}=v^{a_{ij}}E_{jl}K_i, \qquad
K_iF_{jl}=v^{-a_{ij}}F_{jl}K_i,
\\
\label{eq:QB4}
&K'_iE_{jl}=v^{-a_{ij}}E_{jl}K'_i, \qquad
K'_iF_{jl}=v^{a_{ij}}F_{jl}K'_i,
\\
\label{eq:QB5}
&E_{ik}F_{jl}-F_{jl}E_{ik}=\delta_{lk}\delta_{ij}\frac{K_i-K_i'}{v-v^{-1}},
\\
\label{eq:QB7}
&
\sum_{n=0}^{1-a_{ij}}(-1)^n\qbinom{1-a_{ij}}{n} E_{ik}^{1-a_{ij}-n}E_{jl}E_{ik}^n=0,
\\
\label{eq:QB8}
&\sum_{n=0}^{1-a_{ij}}(-1)^n\qbinom{1-a_{ij}}{n} F_{ik}^{1-a_{ij}-n}F_{jl}F_{ik}^n=0,
\text{ for }i\in\Ire, i\neq j.
\end{align}

Note that $K_iK_i'$ are central in $\tUB$ for all $i\in\I$. The quantum generalized Kac-Moody algebra $\UB:=\UB_v(\fg_{A,\bm})$ is generated by the elements $K_i^{\pm1}$, $E_{ik}$, $F_{ik}$ for $i\in\I$, $k=1,\dots,m_i$ subject to \eqref{eq:QBB1}--\eqref{eq:QBB2}, \eqref{eq:QB3}--\eqref{eq:QB8} with $K_i'$ replaced by $K_i^{-1}$.
In particular, the quantum generalized Kac-Moody algebra $\UB$ is isomorphic to the quotient algebra of $\tUB$ modulo the ideal $(K_iK_i'-1\mid i\in\I)$.

Let $\tUB^+$ be the subalgebra of $\tUB$ generated by $E_{ik}$ $(i\in \I,k=1,\dots,m_i)$, $\tUB^0$ be the subalgebra of $\tUB$ generated by $K_i^{\pm1}, (K_i')^{\pm1}$ $(i\in \I)$, and $\tUB^-$ be the subalgebra of $\tUB$ generated by $F_{ik}$ $(i\in \I,k=1,\dots,m_i)$, respectively.
The subalgebras $\UB^+$, $\UB^0$ and $\UB^-$ of $\UB$ are defined similarly. Then both $\tUB$ and $\UB$ have triangular decompositions:
\begin{align}
\label{triandecomp}
\tUB =\tUB^+\otimes \tUB^0\otimes\tUB^-,
\qquad
\UB &=\UB^+\otimes \UB^0\otimes\UB^-.
\end{align}
Clearly, ${\UB}^+\cong\tUB^+$, $\UB^-\cong \tUB^-$, and $\UB^0 \cong \tUB^0/(K_i K_i' -1 \mid   i\in \I)$.

Let $Q=(\I,\Omega)$ be a quiver. Let $\rep_\bfk(Q)$ be the category of finite-dimensional representations of $Q$ over $\bfk$. Then $\rep_\bfk(Q)$ is a hereditary abelian category. We assume that the charge $\bm=(m_i)_{i\in\I}$ satisfies $m_i\leq |\bfk^{g_i}|=q^{g_i}$ for all $i\in\I$.
If $i\in\Ire$ then there exists a unique simple object $S_i\in\rep_\bfk(Q)$ supported at $i$. If $i\in\Iim$ then the set of simple representations supported at $i$ is in bijection with $\bfk^{g_i}$: if $\sigma_1,\dots, \sigma_{g_i}$ denote the simple loops at $i$ then to $\underline{\lambda}=(\lambda_1,\dots,\lambda_{g_i})$ corresponds the simple module $S_i(\underline{\lambda})=(V_j,x_\sigma)$ with $\dim_\bfk V_j=\delta_{ij}$ and $x_{\sigma_l}=\lambda_l$ for $l=1,\dots,g_i$.

Following \cite{KS06}, we choose $\underline{\lambda}_i^{(l)}\in\bfk^{g_i}$ for $l=1,\dots,m_i$ in such a way that
$\underline{\lambda}_i^{(l)}\neq \underline{\lambda}_i^{(l')}$ for $l\neq l'$. We set $S_{il}:=S_{i}(\underline{\lambda}_i^{(l)})$ for $i\in\Iim$ and $l=1,\dots, m_i$; and simply set $S_{i1}:=S_{i}$ for $i\in\Ire$.

In order to realize the quantum generalized Kac-Moody algebras $\tUB$, we modify the definition of semi-derived Ringel-Hall algebras  slightly. The algebra $\cs\widetilde{\ch}(\rep_\bfk (Q))$ is constructed in the same way as $\tMH(\rep_\bfk (Q))$ but with the ideal $I_{\Z_2}$ of $\ch(\cc_{\Z_2}(\rep_\bfk(Q)))$ generated by
\begin{align}
&\{[L]-[K\oplus M]\mid \exists \text{ a short exact sequence }0 \longrightarrow K \longrightarrow L \longrightarrow M \longrightarrow 0 \}
\\
\notag
\bigcup&\{ [K_M]-\prod_{i\in\I}[K_{S_i}]^{a_i},\quad [K^*_M]-\prod_{i\in\I}[K^*_{S_i}]^{a_i}\mid M\in\rep_\bfk(Q) \text{ with }\dimv M=(a_i)_{i\in\I} \}.
\end{align}
In fact, $\cs\widetilde{\ch}(\rep_\bfk (Q))$ is isomorphic to the  quotient algebra of $\tMH(\rep_\bfk (Q))$ modulo the ideal generated by
\begin{align}
\label{eq:ideal2}
[K_M]-\prod_{i\in\I}[K_{S_i}]^{a_i},\quad [K^*_M]-\prod_{i\in\I}[K^*_{S_i}]^{a_i},\text{ for any }M\in\rep_\bfk(Q) \text{ with }\dimv M=(a_i)_{i\in\I}.
\end{align}

\begin{lemma}
\label{lem:basis SH}
We have the following triangular decomposition
\begin{align}
\label{eq:triangdecomp2}
\cs\widetilde{\ch}(\rep_\bfk (Q))=\widetilde{\ch}(\rep_\bfk(Q))\otimes \Q(\sqq)\big[[K_{S_i}],[K^*_{S_i}]\mid i\in\I\big] \otimes \widetilde{\ch}(\rep_\bfk(Q))
\end{align}
\end{lemma}

\begin{proof}
Note that the quotient algebra of the quantum torus $\ct(\rep_\bfk(Q))$ modulo the ideal generated by \eqref{eq:ideal2} is isomorphic to $\Q(\sqq)\big[[K_{S_i}],[K^*_{S_i}]\mid i\in\I\big]$.

Since the Euler form of $\rep_\bfk(Q)$ depends only on dimension vectors of  representations, the proof of Theorem \ref{theorem basis of modified hall algebra} (see the proof of \cite[Theorem 3.20]{LP16}) also works for $\cs\widetilde{\ch}(\rep_\bfk (Q))$ by replacing the quantum torus $\ct(\rep_\bfk(Q))$ by $\Q(\sqq)\big[[K_{S_i}],[K^*_{S_i}]\mid i\in\I\big]$.

Then the same proof of \cite[Theorem 3.25]{LP16} gives \eqref{eq:triangdecomp2}.
\end{proof}

\begin{theorem}
There exists an injective morphism of algebras $\Psi: \tUB\rightarrow \cs\widetilde{\ch}(\rep_\bfk (Q))$ such that
\begin{align}
\Psi(K_i)=[K_{S_i}],\quad \Psi(K_i')=[K_{S_i}^*],\quad \Psi(E_{il})=\frac{1}{q-1}[C_{S_{il}}],\quad \Psi(F_{il})=\frac{-\sqq}{q-1}[C^*_{S_{il}}]
\end{align}
for any $i\in\I$.
\end{theorem}

\begin{proof}
In order to prove the well-definedness of $\Psi$, it is enough to check that $\Psi$ preserves the relations \eqref{eq:QBB1}--\eqref{eq:QBB2}, \eqref{eq:QB3}--\eqref{eq:QB8}. In fact, \eqref{eq:QBB1}--\eqref{eq:QBB2} and \eqref{eq:QB3}--\eqref{eq:QB4} follow from Lemma \ref{lem:KC};
\eqref{eq:QB7}--\eqref{eq:QB8} follow from \cite[Theorem 2.2]{KS06} by using \eqref{eq:Rpm}.

For \eqref{eq:QB5}, if $(i,k)\neq(j,l)$, by definition, we have zero Hom-spaces between $S_{ik}$ and $S_{jl}$. So
$$[C_{S_{ik}}]*[C^*_{S_{jl}}]=[C_{S_{ik}}\oplus C^*_{S_{jl}}]=[C^*_{S_{jl}}]*[C_{S_{ik}}].$$
 If $(i,k)=(j,l)$, a direct computation shows
\begin{align*}
[C_{S_{ik}}]*[C^*_{S_{ik}}]=[C_{S_{ik}}\oplus C^*_{S_{ik}}]+(q-1)[K_{S_i}^*],
\\
[C^*_{S_{ik}}]*[C_{S_{ik}}]=[C_{S_{ik}}\oplus C^*_{S_{ik}}]+(q-1)[K_{S_i}].
\end{align*}
Then \eqref{eq:QB5} follows.

It remains to prove the injectivity of $\Psi$. By \eqref{triandecomp} we have
$\tUB =\tUB^+\otimes \tUB^0\otimes\tUB^-$, correspondingly, by Lemma \ref{lem:basis SH}, we have
\begin{align*}
\cs\widetilde{\ch}(\rep_\bfk (Q))=\widetilde{\ch}(\rep_\bfk(Q))\otimes \Q(\sqq)\big[[K_{S_i}],[K^*_{S_i}]\mid i\in\I\big] \otimes \widetilde{\ch}(\rep_\bfk(Q)).
\end{align*}
Then the injectivity of $\Psi$ follows from \cite[Theorem 2.2]{KS06} by using the same argument of \cite[Theorem 4.9]{Br13}.
\end{proof}

\end{document}